# DFPI, A unified framework for deflated linear solvers: bridging the gap between Krylov subspace methods and Fixed-Point Iterations

Jeremy Kalfoun[a,b,∗], Guillaume Pierrot[b], John Cagnol[a]

[a]*Université Paris-Saclay, CentraleSupélec, MICS Lab & CNRS, Fédération de Mathématiques de CentraleSupélec, Gif-sur-Yvette, 91190, France*
[b]*ICON Technology & Process Consulting, Paris, 75011, France*

**Abstract**

Iterative algorithms are instrumental in modern numerical simulation for solving systems arising from the discretization of PDEs. They face however significant challenges in industrial applications, such as slow convergence, limit cycle oscillations, or iterations blow-up. An ideal preconditioner is rarely available and naive approaches such as Richardson iterations often fail to converge on complex cases, calling for generic sophistications such as deflation techniques and/or Krylov subspaces approaches. However the quest for an optimal general linear solver is still open and a matter of active research. This paper introduces a new theoretical framework, called DFPI (Deflated Fixed Point Iterations) for the iterative solution of linear systems. It unifies several existing acceleration and stabilization techniques such as RPM, BoostConv and Anderson acceleration, and bridges the gap between Richardson iterations and Krylov subspace methods, including GMRES, PCG, BiCGStab and variants. DFPI is structured around two key building blocks : the choice of a projection operator, on the one hand and the trouble vectors recruitment strategy, on the other hand. A general convergence result will be presented, showing the choice of a specific projection operator has minimal impact as long as the projection space remains invariant by the iteration matrix. However when this is not guaranteed, a minimization principle becomes a must have. Finally numerical comparisons will be conducted on a variety of relevant CFD cases.

*Keywords:* Fixed-Point iterations, Deflation methods, Krylov methods, stabilization of solution algorithms, convergence acceleration, CFD
*2000 MSC:* 65B99, 65F10, 65N99

## 1. Introduction

Partial Differential Equations (PDEs) are essential in various modeling contexts, such as Computational Fluid Dynamics (CFD) [1]. As closed-form solutions are generally unavailable, numerical approximation methods like Finite Elements, Finite Differences, or Finite Volumes are used. These methods often lead to large systems of equations, which are then solved iteratively, as direct methods are impractical due to their prohibitive complexity and storage requirements. Over the past decades, the performance of iterative solvers for large systems of the form $Ax = b$ has been a critical area of research.

Two main families of iterative methods exist: fixed-point iterations (Richardson iterations) and Krylov subspace methods. Fixed-point iterations iteratively compute $x^{(k+1)} = Gx^{(k)} + b$ starting from an initial guess $x^{(0)}$. Their convergence depends on the spectral radius of the iteration matrix $G$ [2]. This method is straightforward to implement and has minimal storage requirements, as only $x^{(k)}$ needs to be retained. Krylov subspace methods project the problem onto a sequence of increasingly larger subspaces [2]. The process begins with an initial guess $x$ and the computation of the residual $r = Ax - b$, which is typically non-zero. Krylov spaces are constructed sequentially by iteratively adding vectors formed by applying powers of the matrix $A$ to the residual. This approach might demand significantly

---

more storage than fixed-point iterations depending if a short recurrence is available (CG, BiCG and variants) or not (GMRES and variants) [3, 2, 4]. For SPD matrices, convergence of the (P)CG method is theoretically guaranteed but its convergence rate depends on the condition number of the (preconditioned) system [5] which might lead to residual stagnation in practice. More problematic is that in the general framework, convergence is either not guaranteed at all (BiCG..) or is guaranteed only in case of storage equivalent to the size of $A$ (GMRES) [6] which is unpractical on industrial models, leading to restarted strategies in practice (GMRES(n)) [7, 8] which might in turn faze the convergence.

While relaxation-based methods can offer robust solutions for special matrix types, no universal method exists for solving large, ill-conditioned systems. One common approach to improving convergence involves modifying the problem to limit or regularize troublesome terms, such as non-orthogonal correctors in finite volume methods [9], hard limiting such as velocity damping, or Adjoint Transpose Convection (ATC) damping in adjoint solvers [10, 11]...).

A more satisfactory approach is deflation. It divides the solution space into a "trouble space" where convergence is difficult and an "easygoing space" where it is not, and applies different solving strategies on both of them. The trouble space can either be designed *a priori* (eg. Multigrid, subdomain deflation, Model Order Reduction [12, 13, 14]) or a posteriori via "Krylov-deflation", where dominant eigenvectors are approximated during the solution process. While comprehensive descriptions of deflation exist for Krylov methods [14, 8], no general framework exists for fixed-point iterations. A variety of approaches have however been proposed, such as RPM [15, 16] and BoostConv [17, 18] just to name a few, but each one of them is somewhat lacking generality as it typically does enforce a specific trouble modes recruitment and/or a given projection operator.

In this article, we introduce a new general framework called Deflated Fixed-Point Iterations (DFPI). Focused on the linear context, DFPI unifies several techniques, offering a basis for deriving general convergence results and discussing key aspects like projection operators, trouble space design, and vector recruitment. It also bridges the gap between Krylov methods and fixed-point iterations.

The paper is organized as follows. In section 2, we introduce the deflated fixed-point iterations framework as a combination of Richardson and projection steps, and present an equivalent formulation as a modified preconditioner. In section 3 we present our main convergence result for an invariant deflation subspace, irrespective of the used projector. We then extend this result to quasi-invariant subspaces and to genuinely non-invariant ones under some additional assumption on the projector. Then we address the question of varying projectors due to dynamic trouble spaces, which will enable us to make the connection in between DFPI and Krylov methods. In section 4 we show that this new framework encompasses various well-established solvers, such as RPM or BoostConv. Finally we demonstrate the application of the deflated fixed-point iteration framework in a series of numerical cases.

## 2. Formulation of the DFPI

*2.1. Main formulation*

Let $A \in Gl_N(\mathbb{R})$ be a given non-singular real matrix and $b \in \mathbb{R}^N$ be a known vector. The integer $N$ is the dimension of the system. Typical industrial problems can lead to $N$ having seven to eight digits. Our objective is to find $x^\infty$, the solution to

$$Ax^\infty = b \qquad (1)$$

The *Preconditioned Fixed-Point Iterations method*, also known as the *Richardson method*, is given by

$$\begin{aligned} x^{(n+1)} &= x^{(n)} + P^{-1}(b - Ax^{(n)}) \\ &= x^{(n)} + P^{-1}A(x^\infty - x^{(n)}) \end{aligned} \qquad (2)$$

where the matrix $P \in Gl_N(\mathbb{R})$ (called preconditioner) is selected to ensure that $P^{-1}A$ exhibits favorable spectral properties while remaining computationally efficient to compute and invert. The convergence of the iterations in Eq. (2) depends on the spectral radius of $\text{Id} - P^{-1}A$ [2]. If the spectral radius is strictly less than one, the method converges. However, in practice, finding a preconditioner $P$ that guarantees this condition can be challenging, and some eigenvalues of $\text{Id} - P^{-1}A$ may have a modulus greater than one preventing convergence or close to one, thus slowing it down. Subsequently, we introduce the projection space $\mathbf{Z} = \text{span}\{Z_i\}_{i=1,...,M}$ (here after called "trouble space") where the vectors $Z_i \in \mathbb{R}^N$ are independent and represent problematic modes to be filtered out from the fixed-point iterations. Therefore, we shall operate differently on **Z** and its algebraic complement.

Let $Q_R$ be a projection onto the trouble space **Z**. Any projection operator can be fully defined on the one hand by its image (here **Z**), and on the other hand by its kernel or equivalently the space $\mathbf{Y} = A^{-T}(\ker(Q_R))^\perp$ where $A^{-T}$ is the inverse of the transpose of $A$. Indeed, we have by construction:

$$\forall x \in \mathbb{R}^N, \forall y \in \mathbf{Y}, \left\langle (\mathrm{Id} - Q_R) x, A^T y \right\rangle = 0$$

and we recognize the Petrov-Galerkin projection [2]

$$\forall x \in \mathbb{R}^N, \forall y \in \mathbf{Y}, \langle A Q_R x, y \rangle = \langle A x, y \rangle \tag{3}$$

In this way, the projection of the error vector $x^\infty - x^{(n)}$ is defined based only on the computable residual $A\left(x^\infty - x^{(n)}\right) = b - A x^{(n)}$, which is the sole quantity known in practice. This leads to the operator $Q_R = Z \left(Y^T A Z\right)^{-1} Y^T A$. The subscript "R" refers to a right multiplication by $A$, as opposed to the later introduced operator $Q_L$ for which the multiplication by the matrix $A$ will be performed on the left hand side. The $N \times M$ matrices $Z$ and $Y$ are constructed such that their columns form bases for the subspaces **Z** and **Y**, respectively. Some notable special cases, based on a linear relation between **Z** and **Y** include:

- The *Galerkin projection*, where $\mathbf{Y} = \mathbf{Z}$,
- The *Least Square projection*, where $\mathbf{Y} = A\mathbf{Z}$,
- The *Orthogonal projection*, where $\mathbf{Y} = A^{-T}\mathbf{Z}$

The key idea behind the DFPI method introduced in this work is to solve directly within the trouble space **Z** while iterating on its algebraic complement, referred to as the "easygoing space." This is achieved by projecting the current error $x^\infty - x^{(n)}$ on $\ker(Q_R)$ at the beginning of each iteration. We thus introduce an intermediate step: $x^{(n+\frac{1}{2})} = (\mathrm{Id} - Q_R) x^{(n)} + Q_R x^\infty$. As mentioned earlier, $Q_R x^\infty$ can be calculated as $A x^\infty = b$.

Then: $x^{(n+1)} = x^{(n+\frac{1}{2})} + P^{-1} A (x^\infty - x^{(n+\frac{1}{2})})$. This leads to the first formulation of the Deflated Fixed-Point Iterations (DFPI), whose convergence properties will be discussed in section 3:

$$\begin{cases} x^{(n+\frac{1}{2})} = x^{(n)} + Q_R \left(x^\infty - x^{(n)}\right) \\ x^{(n+1)} = x^{(n+\frac{1}{2})} + P^{-1} A \left(\mathrm{Id} - Q_R\right) \left(x^\infty - x^{(n)}\right) \end{cases} \tag{4}$$

Throughout the paper we will refer to this algorithm as DFPI(**Z**,L). **Z** defines the space we are projecting on, and $L$ the linear operator mapping **Z** into **Y** ($L$ = Id for Galerkin projection, $L = A$ for Least Square Projection...). **Z** can evolve throughout the iterations in some cases, as we will explain later, but $L$ is assume to remain the same. By leveraging the method's degrees of freedom, when it comes to select a specific projection operator and/or the assumed or deduced trouble modes, the approach unifies and generalizes several existing techniques for accelerating and stabilizing fixed-point iterations. Moreover, we will demonstrate that it bridges the gap between Richardson iterations and Krylov subspace methods, as the latter ones can be formulated in the DFPI framework.

### 2.2. Prior and Post projection

It is possible to derive a couple of variants to the baseline algorithm, offering some alternative view on the method. In the linear case, the $Q_R (x^\infty - x^n)$ part does not affect iterations on the algebraic complement $\ker(Q_R)$. In the case where iterations on the algebraic complement converge towards $x^*$, it suffices to correct the $Q_R$ part at the very last iteration, yielding another formulation of the iterations:

$$\begin{cases} x^{(n+1)} = x^{(n)} + P^{-1} A \left(\mathrm{Id} - Q_R\right) \left(x^\infty - x^{(n)}\right) \\ x^\infty = x^* + Q_R (x^\infty - x^*) \end{cases} \tag{5}$$

These two expressions of the DFPI correspond to a projection prior to iteration. In a parallel manner, the projection of the error onto $ker(Q_R)$ can be performed after iterating. This leads to a post-projection formulation:

$$\begin{cases} x^{(n+\frac{1}{2})} = x^{(n)} + P^{-1}A\left(x^\infty - x^{(n)}\right) \\ x^{(n+1)} = x^{(n+\frac{1}{2})} + Q_R\left(x^\infty - x^{(n+\frac{1}{2})}\right) \end{cases} \quad (6)$$

In the linear case, if $x^{(0)}$ satisfies $Q_R x^{(0)} = Q_R x^\infty$, then $\forall n \in \mathbb{N}, Q_R x^{(n)} = Q_R x^\infty$. Consequently, the iterations can be written with a projection at the initialization step only, rather than at the final step.

$$\begin{cases} x^{(0)} = \tilde{x}^{(0)} + Q_R\left(x^\infty - \tilde{x}^{(0)}\right) \\ x^{(n+1)} = x^{(n)} + (\text{Id} - Q_R) P^{-1} A \left(x^\infty - x^{(n)}\right) \end{cases} \quad (7)$$

Permuting fractional and full steps, the pre-projection Eq. (4) and post-projection Eq. (6) formulations prove to be equivalent.

### 2.3. Preconditioner formulation

Let the "left" projector $Q_L$ be defined by $Q_L = A Q_R A^{-1}$ (it is indeed a projector as $Q_L^2 = A Q_R A^{-1} A Q_R A^{-1} = A Q_R Q_R A^{-1} = A Q_R A^{-1} = Q_L$). This projector is of little use in practice due to the presence of $A^{-1}$ but it is of theoretical interest as it allows to reformulate the algorithm as standard Richardson iterations with a modified preconditioner.

**Proposition 1.** *For any projection $Q_R$ such that $A^{-1} Q_L + P^{-1}(\text{Id} - Q_L)$ is non-singular, the DPFI Eq. (4) can be formulated as Fixed-Point iterations using a new preconditioner $\underline{P} = \left(A^{-1} Q_L + P^{-1}(\text{Id} - Q_L)\right)^{-1}$*

PROOF. The main formulation of the DFPI can be written:

$$\begin{aligned}
\underline{x}^{(n+1)} &= \underline{x}^{(n)} + Q_R\left(x^\infty - \underline{x}^{(n)}\right) + P^{-1} A (\text{Id} - Q_R)\left(x^\infty - \underline{x}^{(n)}\right) \\
&= \underline{x}^{(n)} + \left(Q_R + P^{-1} A (\text{Id} - Q_R)\right)\left(x^\infty - \underline{x}^{(n)}\right) \\
&= \underline{x}^{(n)} + \left(A^{-1} Q_L A + P^{-1} (\text{Id} - Q_L) A\right)\left(x^\infty - \underline{x}^{(n)}\right) \\
&= \underline{x}^{(n)} + \left(A^{-1} Q_L + P^{-1} (\text{Id} - Q_L)\right) A \left(x^\infty - \underline{x}^{(n)}\right) \\
&= \underline{x}^{(n)} + \underline{P}^{-1} A \left(x^\infty - \underline{x}^{(n)}\right)
\end{aligned}$$

Conditions for $A^{-1} Q_L + P^{-1}(\text{Id} - Q_L)$ to be non-singular are given section 3.

## 3. Convergence of the DFPI

In this section, we study the convergence of the Deflated Fixed-Point Iteration. Unless explicitly mentioned, no specific assumption on the projector $Q_R$ will be made in the following results.

### 3.1. Main results

In this subsection, we establish convergence results when the subspace $\mathbf{Z}$ is fixed throughout the iterations and invariant with respect to multiplication by $P^{-1}A$. That means that $P^{-1}A\mathbf{Z} \subset \mathbf{Z}$. Subsections 3.2 and 3.3 generalize the results when $\mathbf{Z}$ is not invariant (which is the case for Krylov subspaces). Subsection 3.4 shows that whether a projection on $\mathbf{Z}$ is performed from the start, or only starting at iteration $n \in \mathbb{N}$, from this point, the iterates given by the two variants are equal and the convergence rates are identical thereafter. This means that the convergence results hold with a dynamic projection space, making the link with Krylov methods.

**Lemma 1.** *Let $\mathbf{Z}$ be a subspace of $\mathbb{R}^N$, $Q_R$ be a projection onto $\mathbf{Z}$, $A \in \mathcal{G}l_N(\mathbb{R})$, and $\mathbf{Y} = A^{-T}(\ker(Q_R))^\perp$, then $\forall x \in \mathbb{R}^N, Y^T A (\text{Id} - Q_R) x = 0$.*

PROOF. From Eq. (3), we derive $Y^T A Q_R x = Y^T A x$, thus: $Y^T A (\text{Id} - Q_R) x = 0$. □

**Proposition 2.** *Let $Q_R$ be any projection operator. If the iterations defined by Eq. (4) converge, and the matrix $N = Q_R + P^{-1} A (\text{Id} - Q_R)$ is non-singular, then the iterations converge to the solution of the linear system (1).*

PROOF. In subsection 2.3, the main formulation of the DFPI is expressed as fixed-point iterations using a new preconditioner, denoted as $\underline{P}$. If these iterations converge, they will converge to the solution of Eq. (1) provided that $\underline{P}$ is well-defined and non-singular. This condition is satisfied when $N$ is non-singular. □

**Proposition 3.** *The matrix $N = Q_R + P^{-1} A (\text{Id} - Q_R)$ is non-singular if either of these two conditions is met:*

*(i) $\mathbf{Z}$ is invariant by $P^{-1} A$ (i.e. $P^{-1} A \mathbf{Z} \subset \mathbf{Z}$)*

*(ii) $Y^T P \mathbf{Z}$ is non-singular*

PROOF. Let $x \in \ker(N)$. $Nx = 0 \Rightarrow (\text{Id} - Q_R) N x = 0$.
But $(\text{Id} - Q_R) Q_R = 0$ (from projection's properties). This implies $(\text{Id} - Q_R) P^{-1} A (\text{Id} - Q_R) x = 0$.
So $P^{-1} A (\text{Id} - Q_R) x \in \mathbf{Z}$. Let us denote by $z$ this element.

(i) Assume that $\mathbf{Z}$ is invariant by $P^{-1} A$. By bijectivity of $P^{-1} A$, $P^{-1} A \mathbf{Z} \subset \mathbf{Z} \Rightarrow A^{-1} P \mathbf{Z} \subset \mathbf{Z}$.
Then $(\text{Id} - Q_R) x = A^{-1} P z \in \mathbf{Z}$. Thus $(\text{Id} - Q_R) x = 0$. Substituting in $Nx = 0$, it gives $Q_R x = 0$. And so $x = 0$ and $\underline{P}$ defined and non-singular.

(ii) Assume $Y^T P \mathbf{Z}$ non-singular. $P^{-1} A (\text{Id} - Q_R) x = z \Rightarrow Y^T A (\text{Id} - Q_R) x = Y^T P z$. There exists $C \in \mathbb{R}^M$ such that $z = ZC$ and from lemma 1, $Y^T A (\text{Id} - Q_R) x = 0$. Then $Y^T P Z C = 0$, and as $Y^T P \mathbf{Z}$ is non-singular, it yields $z = 0$.
Again, substituting in $Nx = 0$ implies $Q_R x = 0$, so $x = 0$ and $\underline{P}$ defined and non-singular. □

These results state conditions under which the limit of DFPI iterations, when it exists, matches the expected linear system solution. Next, we will derive conditions for the iterations to actually converge.

**Remark 1.** With a similar condition than the one stated for convergence with the preconditioner $P$, iterations using the preconditioner $\underline{P}$, will converge if and only if the following condition on the spectral radius holds:

$$\forall \lambda \in \text{Sp}(\text{Id} - (Q_R + P^{-1} A (\text{Id} - Q_R))), |\lambda| < 1$$

Let us consider a space $\mathbf{Z}$, that contains all the generalized eigenvectors associated with the troublesome eigenvalues of $\text{Sp}(\text{Id} - P^{-1} A)$, preventing or disturbing convergence of the iterations with preconditioner $P$. Denote these eigenvalues as $\lambda_i$ for $1 \leq i \leq M$, where $\text{Sp}(\text{Id} - P^{-1} A) = \{\lambda_i\}_{1 \leq i \leq N}$ and $M < N$. It means that $\forall i \in [\![M+1, N]\!], |\lambda_i| < 1$.
We state the following theorem for any projection:

**Theorem 1.** *Let $\mathbf{Z}$ be a $(P^{-1} A)$-invariant subspace (i.e. $P^{-1} A \mathbf{Z} \subset \mathbf{Z}$) and let $Q_R$ be a projection onto $\mathbf{Z}$. Then the spectrum of the matrix $\text{Id} - (Q_R + P^{-1} A (\text{Id} - Q_R))$ is given by:*

$$\text{Sp}(\text{Id} - (Q_R + P^{-1} A (\text{Id} - Q_R))) = (0, 0, ..., 0, \lambda_{M+1}, ..., \lambda_N)$$

*where $\lambda_{M+1}, \ldots, \lambda_N$ are the non-deflated eigenvalues of $P^{-1} A$. Consequently, iterations with the preconditioner $\underline{P}$ converge at a rate determined by the largest non-deflated eigenvalue in modulus.*

PROOF. Consider the direct sum decomposition $\mathbb{R}^N = \text{Im}(Q_R) \oplus \ker(Q_R)$. Let $\mathcal{B}'$ be a basis for $\mathbb{R}^N$ constructed from bases of $\text{Im}(Q_R)$ and $\ker(Q_R)$. In this basis, the matrices $M = P^{-1} A$ and $N = Q_R + P^{-1} A (\text{Id} - Q_R)$ can be written as:

Without loss of generality on the eigenvalues, matrices are from here written in this basis. Let's consider the two matrices $M = P^{-1} A$ and $N = Q_R + P^{-1} A (\text{Id} - Q_R)$. As $\mathbf{Z}$ is invariant by $P^{-1} A$, $(\text{Id} - Q_R) P^{-1} A Q_R = 0$ and:

$$M = \begin{pmatrix} (Q_R P^{-1} A Q_R)|_{\text{Im}(Q_R)} & (Q_R P^{-1} A (\text{Id} - Q_R))|_{\ker(Q_R)} \\ 0 & ((\text{Id} - Q_R) P^{-1} A (\text{Id} - Q_R))|_{\ker(Q_R)} \end{pmatrix},$$

$$N = \begin{pmatrix} Q_R|_{\text{Im}(Q_R)} & (Q_R P^{-1} A (\text{Id} - Q_R))|_{\ker(Q_R)} \\ 0 & ((\text{Id} - Q_R) P^{-1} A (\text{Id} - Q_R))|_{\ker(Q_R)} \end{pmatrix}$$

Since **Z** is invariant under $P^{-1}A$, we have $(\text{Id} - Q_R)P^{-1}AQ_R = 0$. Thus, the eigenvalues of $M$ and $N$ are determined by the eigenvalues of their respective diagonal blocks.

The eigenvalues of the bottom right block of $M$ and $N$ are identical, while the eigenvalues of $\left(Q_R P^{-1} A Q_R\right)|_{\text{Im}(Q_R)}$ on $M$ correspond to the eigenvalues of $Q_R|_{\text{Im}(Q_R)}$ on $N$, which are all ones.

Therefore, the spectrum of $\text{Id} - (Q_R + P^{-1}A(\text{Id} - Q_R))$ is:

$$\text{Sp}\left(\text{Id} - \left(Q_R + P^{-1} A \left(\text{Id} - Q_R\right)\right)\right) = (0, 0, ..., 0, \lambda_{M+1}, ..., \lambda_N)$$

where $\lambda_{M+1}, \ldots, \lambda_N$ are the eigenvalues of $\left((\text{Id} - Q_R) P^{-1} A (\text{Id} - Q_R)\right)|_{\ker(Q_R)}$. By hypothesis, **Z** is chosen such that $\forall i \in [\![M+1, N]\!]$, $|\lambda_i| < 1$ ensuring convergence of iterations with the preconditioner $\underline{P}$. □

**Remark 2.** Theorem 1 extends the results presented in [15] in the linear framework by generalizing the theorem to any projection.

Moreover, we have an additional result on the structure of the deflated eigenvectors, that is the generalized eigenvectors of the matrix $N = Q_R + P^{-1}A(\text{Id} - Q_R)$ as compared to $M = P^{-1}A$.

Let us introduce the Jordan chain of the matrix $M$ associated with the $i$-th eigenvalue $\lambda_i$. An eigenvector corresponding to $\lambda_i$ denoted $\epsilon_{i0}$ is the first element of the chain, then come the vectors $\epsilon_{ij}$ for all $j \in \{1, m_i - 1\}$, if the following relation holds: $M\epsilon_{ij} = \lambda_i \epsilon_{ij} + \epsilon_{ij-1}$. The chain is denoted

$$C(M, \lambda_i, (\epsilon_{ij})_{1 \leq j \leq m_i - 1}) \in \mathbb{R}^{m_i - 1}$$

Its elements are the generalized eigenvectors.

**Theorem 2.** *Under the hypothesis of Theorem 1, there exists a bijection between the generalized eigenvectors of the matrices $M$ and $N$. Specifically, the generalized eigenvectors $(e_{ij})_{1 \leq i \leq n,\ 0 \leq j \leq m_i - 1}$ of $N$ can be expressed as*

$$e_{ij} = \mu_i (\text{Id} - Q_R) \epsilon_{ij} + z_{ij}$$

*where $(\epsilon_{ij})_{1 \leq i \leq n,\ 0 \leq j \leq m_i - 1}$ are the generalized eigenvectors of $M$, $\mu_i \in \mathbb{R}$, and $z_{ij} \in \mathbf{Z}$. In particular, when $\epsilon_{ij} \in \mathbf{Z}$, we have $e_{ij} = \mu_i \epsilon_{ij}$, which is associated with the eigenvalue 1.*

The proof is given in Appendix Appendix A.

### 3.2. Robustness of the results

In practice **Z** is rarely fully invariant under the action of $P^{-1}A$ due to finite precision in the approximation of eigenvectors. As a consequence, the matrix $N = Q_R + P^{-1}A(\text{Id} - Q_R)$ does not exhibit an exact triangular structure.

In the following we denote by **Z** the invariant space with respect to $P^{-1}A$, consisting of $m$ eigenvectors $(z_1, ..., z_m)$, and by $\mathbf{Z}_u = (z_1 + \delta_1, ..., z_m + \delta_m)$ a slight perturbation of it (where $\delta_i \in \mathbb{C}^n$). We then express $N$ as $N = N_T + \delta N$, where $\delta N \in \mathcal{M}_n(\mathbb{C})$ is the perturbation and $N_T$ is the triangular matrix:

$$N_T = \begin{pmatrix} Q_R|_{\text{Im}(Q_R)} & B \\ 0 & \left((\text{Id} - Q_R) P^{-1} A (\text{Id} - Q_R)\right)|_{\ker(Q_R)} \end{pmatrix}$$

with $Q_R$ the projection on **Z**. The spectrum of $N_T$ is given by Theorem 1. For $X \in Gl_n(\mathbb{C})$ let's denote $\kappa_p = \|X\|_p \|X^{-1}\|_p$ the condition number of $X$ associated with the matrix $p$-norm $\|.\|_p$. Let's consider $\lambda$ the spectral radius of $N_T$, corresponding to the largest eigenvalue in modulus associated with an eigenvector that is not in **Z** (Theorem 1). And finally, let's write $V$ the matrix of eigenvectors of $N_T$ and $\rho(N)$ the spectral radius of $N$. Then, we have the following result.

**Proposition 4.** $\rho(N) \leq \lambda + \kappa_p(V) \|\delta N\|_p$

Proof. The Bauer-Fike theorem [19] ensures that $\forall \mu \in \mathrm{Sp}(N), \exists \lambda_i \in \mathrm{Sp}(N_T)$ such that $|\mu - \lambda_i| \leq \kappa_p(V)\|\delta N\|_p$. Hence $||\mu| - |\lambda_i|| \leq |\mu - \lambda_i| \leq \kappa_p(V)\|\delta N\|_p$. If $|\mu| \leq |\lambda_i|$, then $|\mu| \leq \lambda$. Else, $|\mu| - |\lambda_i| \leq \kappa_p(V)\|\delta N\|_p$ by hypothesis. It follows $|\mu| \leq |\lambda_i| + \kappa_p(V)\|\delta N\|_p \leq \lambda + \kappa_p(V)\|\delta N\|_p$. Hence the result. □

In other words, the distance of the eigenvalues of the perturbed matrix to the spectrum of $N_T$ induced by the loss of the stability property is bounded by the norm of the perturbation times a coefficient. If the norm of the perturbation is sufficiently small compared to the distance of the spectrum of $N_T$ to the unit circle, eigenvalues of $N$ will lay inside the unit circle. Let us suppose that eigenvectors in **Z** have been chosen such that $\lambda < 1 - \delta$ with $\delta > 0$.

**Corollary 1.** *If $\delta > \kappa_p(V)\|\delta N\|_p$, then iterations associated to $N$ converge. The spectral radius governing the convergence rate is smaller than $1 - \delta + \kappa_p(V)\|\delta N\|_p$.*

In practice $\delta$ can be made small enough by adjusting the tolerance when recruiting approximate eigenvectors.

*3.3. Non-invariant space and minimization*

Let's now take a different perspective and consider the case where **Z** is genuinely not invariant by $P^{-1}A$ (eg, is a Krylov subspace) but contains an invariant subspace. Let $S \in \mathcal{M}_N(\mathbb{R})$ be a symmetric positive definite matrix. $S$ defines an inner product on $\mathbb{R}^N$: $\forall x, y \in \mathbb{R}^N, \langle x, y\rangle_S = y^T S x$.

**Proposition 5.** *Let $L$ be the linear mapping associated to $S$ defined by $L = A^{-T}S$. Let $\mathbf{I} \subset \mathbf{Z}$ be the largest contained invariant subspace by $P^{-1}A$. DFPI(**Z**,L) converge at a higher rate than DFPI(**I**,L) for the norm induced by the inner product associated to $S$.*

Proof. The projections induced by $L$ on **Z** and **I** are minimizing the norm induced by $\langle . , . \rangle$. As $\mathbf{I} \subset \mathbf{Z}$, the norm of the error produced by DFPI(**Z**,L) is necessarily smaller than the norm of the error given by DFPI(**I**,L). □

**Corollary 2.** *DFPI(**Z**,L) converge at a higher rate than the one given by the largest eigenvalue in modulus associated to an eigenvector not in **I**.*

Proof. It follows directly from Theorem 1. □

This result gives insight on why GMRES and CG (that use orthogonal projections based on $A^T A$ and $A$, respectively) converge seamlessly while BiCG, that is lacking a minimization principle, exhibits an erratic convergence behaviour. This shall be discussed in more details subsequently.

*3.4. Vectors selection and dynamic projection*

In the previous subsections, we presented DFPI iterations based on a given, static trouble space, which corresponds to the situation where it is defined upfront, in an *a priori* manner (for example, based on coarse-grid agglomeration). In practice, however, when we want to use approximate eigenvectors, in order to leverage the result given by Theorem 1, they are to be computed on the fly through generalized power iterations (using vectors increments $x^{(k+1)} - x^{(k)}$ as proxies). Different eigenvectors are potentially captured at different iterations, until we reach either a maximal storage size or a small enough eigenvalue. In the process, two options are possible: one can either modify the projection with the enriched trouble space at each iteration or run unaltered Richardson iterations until some satisfactory projection space has been assembled.

Let's consider the resulting two sequences $(u^{(i)})_{i \in \mathbb{N}}$ and $(v^{(i)})_{i \in \mathbb{N}}$ such that $u^{(0)} = v^{(0)} = x^{(0)}$, $n \in \mathbb{N}$ and

$$\forall k < n, \quad \begin{cases} v^{(k+1)} = v^{(k)} + P^{-1}A\left(x^\infty - v^{(k)}\right) \\ u^{(k+1)} = u^{(k)} + Q_u^{(k)}\left(x^\infty - u^{(k)}\right) + P^{-1}A\left(\mathrm{Id} - Q_u^{(k)}\right)\left(x^\infty - u^{(k)}\right) \end{cases} \quad (8)$$

$$\forall k \geq n, \quad \begin{cases} v^{(k+1)} = v^{(k)} + Q_v^{(n)}\left(x^\infty - v^{(k)}\right) + P^{-1}A\left(\mathrm{Id} - Q_v^{(n)}\right)\left(x^\infty - v^{(k)}\right) \\ u^{(k+1)} = u^{(k)} + Q_u^{(n)}\left(x^\infty - u^{(k)}\right) + P^{-1}A\left(\mathrm{Id} - Q_u^{(n)}\right)\left(x^\infty - u^{(k)}\right) \end{cases} \quad (9)$$

where $Q_u^{(k)}$ is a projection on $\mathbf{Z}_u^{(k)} = \text{span}(u^{(i)} - u^{(i-1)})_{1 \leq i \leq k}$ for $k \leq n$ and $Q_v^{(n)}$ is a projection on $\mathbf{Z}_v^{(n)} = \text{span}(v^{(i)} - v^{(i-1)})_{1 \leq i \leq n}$, based on the same linear operator $L$. Finally, let

$$\mathbb{K}^n = \left\{ \tilde{r}_0, P^{-1}A\tilde{r}_0, \left(P^{-1}A\right)^2 \tilde{r}_0, ..., \left(P^{-1}A\right)^{n-1} \tilde{r}_0 \right\}$$

be the Krylov subspace constructed from $\tilde{r}_0 = P^{-1}A\left(x^\infty - x^0\right)$

**Proposition 6.** $\forall n \in \mathbb{N}, \mathbf{Z}_u^{(n)} = \mathbf{Z}_v^{(n)} = \mathbb{K}^n$.

PROOF. From Eq. (8), the error at iteration $n$ is: $x^\infty - v^{(n)} = \left(\text{Id} - P^{-1}A\right)\left(x^\infty - v^{(n-1)}\right)$, which lead to the increment:

$$v^{(n)} - v^{(n-1)} = v^{(n-1)} - v^{(n-2)} - P^{-1}A(v^{(n-1)} - v^{(n-2)}).$$

Assume $\forall k < n, \mathbf{Z}_v^{(k)} \subset \mathbb{K}^k$, then: $v^{(n-1)} - v^{(n-2)} \in \mathbf{Z}_v^{(n-1)} \subset \mathbb{K}^{n-1}$ and $P^{-1}A\left(v^{(n-1)} - v^{(n-2)}\right) \in \mathbb{K}^n$.
We have $\mathbf{Z}_v^{(n-1)} \subset \mathbb{K}^{n-1} \subset \mathbb{K}^n$ and the last vector $v^{(n)} - v^{(n-1)} \in \mathbb{K}^n$, then $\mathbf{Z}_v^{(n)} \subset \mathbb{K}^n$.

Now assume $\forall k < n$, $\mathbb{K}^{n-1} \subset \mathbf{Z}_v^{(n-1)}$ (and so $\mathbb{K}^k = \mathbf{Z}_v^{(k)}$). As $\left(P^{-1}A\right)^{n-2} \tilde{r}_0 \in \mathbb{K}^{n-1} \subset \mathbf{Z}_v^{n-1}$ it can be decomposed in $\mathbf{Z}_v^{(n-1)}$ using coefficients $\mu_k$, $1 \leq k \leq n-1$. We have

$$\left(P^{-1}A\right)^{n-1} \tilde{r}_0 = \mu_{n-1} P^{-1}A\left(v^{(n-1)} - v^{(n-2)}\right) + \sum_{k=1}^{n-2} \mu_k P^{-1}A\left(v^{(k)} - v^{(k-1)}\right).$$

All the terms inside the summation belong to $\mathbb{K}^{n-1} \subset \mathbf{Z}_v^{(n-1)}$. From the increment equation:

$$(P^{-1}A)^{n-1}\tilde{r}_0 = \mu_{n-1}(v^{(n-1)} - v^{(n-2)}) - \mu_{n-1}(v^{(n)} - v^{(n-1)}) + y^{(n-1)}$$

where $y^{(n-1)} \in \mathbb{K}^{n-1} = \mathbf{Z}_v^{(n-1)}$. So $\left(P^{-1}A\right)^{n-1} \tilde{r}_0 \in \mathbf{Z}_v^{(n-1)}$.
  Initialization is trivial. So finally $\forall n \in \mathbb{N}, \mathbf{Z}_v^{(n)} = \mathbb{K}^n$.
  Then, one can notice that $\forall k \leq n, u^{(k)} - u^{(k-1)}$ is just a shift of $v^{(k)} - v^{(k-1)}$ in the Krylov subspace $\mathbb{K}^k$. This can be shown by recurrence. Initialization is trivial, and let suppose that $\forall i < k, u^{(i)} - u^{(i-1)} = v^{(i)} - v^{(i-1)} + y^{(i)}$ where $y^{(i)} \in \mathbb{K}^i$. From Eq. (8), the increment is written:

$$u^{(k)} - u^{(k-1)} = u^{(k-1)} - u^{(k-2)} + Q_u^{(k-1)}\left(x^\infty - u^{(k-1)}\right) - Q_u^{(k-2)}\left(x^\infty - x^{(k-2)}\right)$$
$$- P^{-1}A\left(u^{(k-1)} - u^{(k-2)} + Q_u^{(k-1)}\left(x^\infty - u^{(k-1)}\right) - Q_u^{(k-2)}\left(x^\infty - u^{(k-2)}\right)\right).$$

From the assumption, $u^{(k-1)} - u^{(k-2)} - P^{-1}A\left(u^{(k-1)} - u^{(k-2)}\right) = v^{(k-1)} - v^{(k-2)} + y^{(k-1)} - P^{-1}A\left(v^{(k-1)} - v^{(k-2)} + y^{(k-1)}\right)$.
The increment for $v$ gives $v^{(k)} - v^{(k-1)} = v^{(k-1)} - v^{(k-2)} - P^{-1}A(v^{(k-1)} - v^{(k-2)})$. And so:

$$u^{(k-1)} - u^{(k-2)} - P^{-1}A\left(u^{(k-1)} - u^{(k-2)}\right) = v^{(k)} - v^{(k-1)} + y^{(k-1)} - P^{-1}Ay^{(k-1)}.$$

All the other terms belong to $\mathbb{K}^k$, and so $u^{(k)} - u^{(k-1)} = v^{(k)} - v^{(k-1)} + y^{(k)}$ where $y^{(k)} \in \mathbb{K}^k$. Hence the equality of the three subspaces. □

Now let's consider the iteration $n$ from which $(v^{(k)})_{k \in \mathbb{N}}$ starts being projected and let's consider the fractional-iterate:

$$\begin{cases} v^{(n+1/2)} = \left(\text{Id} - Q_v^{(n)}\right) v^{(n)} + Q_v^{(n)} x^\infty \\ u^{(n+1/2)} = \left(\text{Id} - Q_u^{(n)}\right) u^{(n)} + Q_u^{(n)} x^\infty \end{cases} \tag{10}$$

Let's also consider the sequence given by Eq. (4), and $x^{(1/2)} = \left(\text{Id} - Q_x^{(n)}\right) x^{(0)} + Q_x^{(n)} x^\infty$ where $Q_x^{(n)}$ is the projection on $\mathbb{K}^n$ based on the same linear operator $L$ as $Q_u^{(n)}$ and $Q_v^{(n)}$.

**Proposition 7.** *For all $n \in \mathbb{N}$, we have $v^{(n+1/2)} = u^{(n+1/2)} = x^{(1/2)}$.*

Proof. From Prop. 6, $\exists y^{(n)} \in \mathbb{K}^n$, $v^{(n)} - v^{(n-1)} = y^{(n)}$. Substrating by $x^{(0)}$ and projecting:

$$\left(\mathrm{Id} - Q_v^{(n)}\right)\left(v^{(n)} - x^{(0)}\right) = \left(\mathrm{Id} - Q_v^{(n)}\right)\left(v^{(n-1)} - x^{(0)}\right).$$

And so by recurrence, as $\left(\mathrm{Id} - Q_v^{(n)}\right)\left(v^{(1)} - x^{(0)}\right) = \left(\mathrm{Id} - Q_v^{(n)}\right)\tilde{r}_0 = 0$ (as $\tilde{r}_0 \in \mathbb{K}^n = \mathbf{Z}_u^{(n)}$), we have

$$v^{(n+1/2)} = \left(\mathrm{Id} - Q_v^{(n)}\right) x^{(0)} + Q_v^{(n)} x^{\infty}$$

which is the same vector as $x^{(1/2)}$. The same reasoning can be applied to $u^{(n)}$. □

This means that steadily projecting on a growing troublespace from iteration 1 to $n$ (as do Krylov methods); performing $n$ vanilla Richardson iterations and then deflating at the last step on the troublespace gathered so far; or doing only one iteration of DFPI with a projection on this very space, gives equivalent iterates between the three methods at all subsequent steps ($u^{(n+1)} = v^{(n+1)} = x^{(1)}$ and so on...), keeping the troublespace constant from this step on. In this case, the convergence rate afterward is upper bounded by Prop. 5. On the other hand, if one keeps increasing the troublespace afterwards, such as in Krylov methods, and some minimization principle is available (as in GMRES or CG), the convergence rate remains at least as good (based on a similar idea than Prop. 5's proof). The key insight is that any iterate of a given Krylov method can always been matched by the outcome of the DFPI algorithm using the method's projection and the gathered Krylov sequence as the trouble space. This allows us to study the different methods in the following sections regardless of the projection space being static or dynamic and makes the link between DFPI and Krylov methods.

## 4. Incorporation of existing methods

In this section, we show how some well-known methods fit within the Deflated Fixed Point Iterations framework. Although some of them are defined in the more general, non-linear case, we will restrict our analysis to the linear setting. Furthermore, on the basis of former section comments, one shall assume that the projection space is static and defer discussions regarding the recruitment of trouble modes to the next paragraph.

### 4.1. RPM in the DFPI framework

The Recursive Projection Method (or RPM) was first introduced by Schroff and Keller [15] to stabilize fixed-point iterative methods by using Newton iterations on the troublespace and fixed-point iterations on the algebraic complement. The technique is described and evaluated in [16]. The method is based on the fixed-point iterations, which, in the case of Preconditioned Richardson Iterations, is written:

$$x^{(n+1)} = F(x^{(n)}) = x^{(n)} + P^{-1}A\left(x^{\infty} - x^{(n)}\right) \tag{11}$$

By introducing the subspace **Z** formed by the eigenvectors associated to dominant eigenvalues, and the orthogonal projection $O$ on **Z**, the iterative scheme is modified as follows in the usual case of additive RPM:

$$\begin{cases} x^{(n)} = \hat{x}^{(n)} + \tilde{x}^{(n)} \\ \hat{x}^{(n+1)} = OF(\hat{x}^{(n+1)} + \tilde{x}^{(n)}) \\ \tilde{x}^{(n+1)} = (\mathrm{Id} - O) F(\hat{x}^{(n)} + \tilde{x}^{(n)}) \end{cases} \tag{12}$$

In our work, we also present multiplicative RPM, where $\hat{x}^{(n+1)}$ is used to calculate $\tilde{x}^{(n+1)}$, as it has been calculated in the previous step:

$$\tilde{x}^{(n+1)} = (\mathrm{Id} - O) F(\hat{x}^{(n+1)} + \tilde{x}^{(n)})$$

Correspondence between DFPI and multiplicative RPM is more straightforward as it does not require the introduction of a new projection $Q_L$ and the hypothesis that **Z** is invariant by $P^{-1}A$ (inherent to RPM in any case).

In classical non-linear RPM, $\hat{x}^{(n+1)}$ is calculated implicitly using a Newton iteration. In the linear case, only one Newton step is necessary.

Developing the fixed-point iteration for $\hat{x}^{(n+1)}$:

$$\hat{x}^{(n+1)} = O\left(\hat{x}^{(n+1)} + \tilde{x}^{(n)} + P^{-1}A\left(x^\infty - \tilde{x}^{(n)} - \hat{x}^{(n+1)}\right)\right)$$

Using that $O\hat{x}^{(n+1)} = \hat{x}^{(n+1)}$, $O\tilde{x}^{(n)} = 0$ and $\tilde{x}^{(n)} = x^{(n)} - \hat{x}^{(n)}$ leads to

$$O\left[P^{-1}A\left(x^\infty - x^{(n)}\right) - P^{-1}A\left(\hat{x}^{(n+1)} - \hat{x}^{(n)}\right)\right] = 0$$

It means that $\forall j \in [\![1, m]\!]$, $\left(P^{-1}A\left[\left(x^\infty - x^{(n)}\right) - \left(\hat{x}^{(n+1)} - \hat{x}^{(n)}\right)\right], Z_j\right) = 0$ with $\hat{x}^{(n+1)} - \hat{x}^{(n)} \in \mathbf{Z}$ which means that $\hat{x}^{(n+1)} - \hat{x}^{(n)}$ is the Galerkin projection of $x^\infty - x^{(n)}$ associated to $P^{-1}A$ on $\mathbf{Z}$.

$$\hat{x}^{(n+1)} = \hat{x}^{(n)} + Q_\xi\left(x^\infty - x^{(n)}\right) \tag{13}$$

In the case of additive RPM, $\tilde{x}^{(n+1)} = (\mathrm{Id} - O) F(\hat{x}^{(n)} + \tilde{x}^{(n)})$, the left projection $Q_L$ such that $Q_L P^{-1} A = P^{-1} A Q_\xi$ is introduced.

$$\tilde{x}^{(n+1)} = \tilde{x}^{(n)} + (\mathrm{Id} - O) Q_L P^{-1} A \left(x^\infty - x^{(n)}\right) + (\mathrm{Id} - O)(\mathrm{Id} - Q_L) P^{-1} A \left(x^\infty - x^{(n)}\right)$$

By definition $Q_L$ is a projection onto $P^{-1}A\mathbf{Z}$. Under the assumption that $\mathbf{Z}$ is invariant by $P^{-1}A$, we have

$$(\mathrm{Id} - O) Q_L P^{-1} A \left(x^\infty - x^{(n)}\right) = 0$$
$$\tilde{x}^{(n+1)} = \tilde{x}^{(n)} + (\mathrm{Id} - O) P^{-1} A \left(\mathrm{Id} - Q_\xi\right)\left(x^\infty - x^{(n)}\right)$$

But from 1 and the definition of the Galerkin projection $Y = Z$,

$$OP^{-1}A\left(\mathrm{Id} - Q_\xi\right)\left(x^\infty - x^{(n)}\right) = Z\left(Z^T Z\right)^{-1}\left[Z^T P^{-1} A \left(\mathrm{Id} - Q_\xi\right)\left(x^\infty - x^{(n)}\right)\right] = 0$$

Therefore $\tilde{x}^{(n+1)} = \tilde{x}^{(n)} + P^{-1}A\left(\mathrm{Id} - Q_\xi\right)\left(x^\infty - x^{(n)}\right)$. Taking

$$\begin{cases} x^{(n+\frac{1}{2})} = \tilde{x}^{(n)} + \hat{x}^{(n+1)} \\ x^{(n+1)} = \tilde{x}^{(n+1)} + \hat{x}^{(n+1)} \end{cases}$$

$$\begin{cases} x^{(n+\frac{1}{2})} = x^{(n)} + Q_\xi\left(x^\infty - x^{(n)}\right) \\ x^{(n+1)} = x^{(n+\frac{1}{2})} + P^{-1}A\left(\mathrm{Id} - Q_\xi\right)\left(x^\infty - x^{(n)}\right) \end{cases}$$

which correspond to the main formulation of the DFPI.

In the case of multiplicative RPM, $\tilde{x}^{(n+1)} = (\mathrm{Id} - O) F(\hat{x}^{(n+1)} + \tilde{x}^{(n)})$. Therefore

$$\tilde{x}^{(n+1)} = \tilde{x}^{(n)} + (\mathrm{Id} - O) P^{-1} A \left(x^\infty - \tilde{x}^{(n)} - \hat{x}^{(n+1)}\right)$$
$$= \tilde{x}^{(n)} + (\mathrm{Id} - O) P^{-1} A \left(x^\infty - \tilde{x}^{(n)} - \hat{x}^{(n)} - Q_\xi\left(x^\infty - x^{(n)}\right)\right)$$
$$= \tilde{x}^{(n)} + (\mathrm{Id} - O) P^{-1} A \left(\mathrm{Id} - Q_\xi\right)\left(x^\infty - x^{(n)}\right)$$

This results in the same conclusion as for the additive case without the need for using $Q_L$ and stability hypothesis.

*4.2. BoostConv in the DFPI framework*

The BoostConv method, introduced by Citro et al. [17], minimizes the residual norm at each step and is inspired by Krylov methods to stabilize the computation of unstable steady states. Its implementation is further discussed in [18], which reports satisfactory results. The method is detailed in [17] and starts with an iterative algorithm to solve Eq. (1):

$$x^{(n+1)} = x^{(n)} + P^{-1} r^{(n)}$$

with $r^{(n)} = b - Ax^{(n)}$. This corresponds to the Preconditioned Fixed-Point Iterations. These iterations, and more precisely the residual are modified by the method as follow:

$$x^{(n+1)} = x^{(n)} + P^{-1}\xi_1^{(n)} \tag{14a}$$
$$\xi_1^{(n)} = \xi_0^{(n)} + \rho^{(n)} \tag{14b}$$
$$\rho^{(n)} = r^{(n)} - AP^{-1}\xi_0^{(n)} \tag{14c}$$
$$\xi_0^{(n)} = \mathrm{argmin}_{\xi \in \mathbf{Z}_0}\|r^{(n)} - AP^{-1}\xi\|_2^2 \tag{14d}$$

where $\mathbf{Z}_0$ is spanned by a certain set of vectors ($\mathbf{u}_i$ in [17]).

The minimization problem defined in Eq. (14d) can be written differently, with a change of variable $\zeta = P^{-1}\xi$ and writing $r^{(n)} = A(x^\infty - x^{(n)})$:

$$P^{-1}\xi_0^{(n)} = \mathrm{argmin}_{\zeta \in \mathbf{Z}}\|A\left((x^\infty - x^{(n)}) - \zeta\right)\|_2^2 \tag{15}$$

with $\mathbf{Z}$ being the subspace given by $P^{-1}\mathbf{Z}_0$. The solution of this minimization problem is nothing but the Least Square Projection of $(x^\infty - x^{(n)})$ on $\mathbf{Z}$. Hence, the solution is written:

$$P^{-1}\xi_0^{(n)} = Q_R^{LSQ}(x^\infty - x^{(n)}) \tag{16}$$

with the Least square projection operator on the subspace $\mathbf{Z}$:

$$Q_R^{LSQ} = Z\left(Z^T A^T A Z\right)^{-1} Z^T A^T A$$

Then the correction $\rho^{(n)} = r^{(n)} - AP^{-1}\xi_0^{(n)}$ is applied to $\xi_0^{(n)}$ such that

$$\xi_1^{(n)} = \xi_0^{(n)} + r^{(n)} - AP^{-1}\xi_0^{(n)}$$
$$= PQ_R(x^\infty - x^{(n)}) + A(x^\infty - x^{(n)}) - AQ_R(x^\infty - x^{(n)})$$

(with $Q_R = Q_R^{LSQ}$) Hence, the modified iteration in Eq. (14a) can be read as:

$$x^{(n+1)} = x^{(n)} + P^{-1}\xi_1^{(n)}$$
$$= x^{(n)} + Q_R(x^\infty - x^{(n)}) + P^{-1}A\,(\mathrm{Id} - Q_R)(x^\infty - x^{(n)})$$

Finally, the BoostConv iterations can be written as follow:

$$\begin{cases} x^{(n+\frac{1}{2})} = x^{(n)} + Q_R\left(x^\infty - x^{(n)}\right) \\ x^{(n+1)} = x^{(n+\frac{1}{2})} + P^{-1}A\,(\mathrm{Id} - Q_R)(x^\infty - x^{(n)}) \end{cases} \tag{17}$$

which correspond exactly to the iterations presented in Eq. (4).

Boostconv is very similar in essence to GMRES, except that the inner product in use is based on $A^T A$ and not $(P^{-1}A)^T(P^{-1}A)$ (furthermore as we will see later on, the projection space corresponds at each iteration to the same Krylov subspace than GMRES)

### 4.3. GMRES in the DFPI framework

The Generalized Minimal Residual (GMRES) method has been introduced by Saad and Schultz [4]. This is a Krylov method that minimizes the residual over a Krylov subspace at each iteration (without restart). A thorough description of GMRES coupled with preconditioning techniques is presented in [2]. GMRES is used on various domains to solve linear systems that might not be stable using fixed point iterations [20, 21, 22].

The method, with left preconditioning, to solve the linear system (1) is based on the Krylov subspaces:

$$\forall n \in \mathbb{N},\ \mathbb{K}^n = \left\{\tilde{r}_0, P^{-1}A\tilde{r}_0, \left(P^{-1}A\right)^2 \tilde{r}_0, ..., \left(P^{-1}A\right)^{n-1} \tilde{r}_0\right\}$$

where $\tilde{r}_0 = P^{-1}r_0 = P^{-1}A\left(x^\infty - x^{(0)}\right)$ with $x^{(0)}$ the initial guess.

Arnoldi iterations are used to construct an orthonormal basis of $\mathbb{K}^n$ during the process through the introduction of an upper Hessenberg matrix. But irrespective of these technicalities, at iteration $k \in \mathbb{N}$ GMRES solves the following minimization problem:

$$z^{(k)} = \mathrm{argmin}_{z \in \mathbb{K}^k} \|P^{-1}b - P^{-1}A\left(x^{(0)} + z\right)\|_2^2$$

GMRES iterate $x^{(k)}$ is given by $x^{(k)} = x^{(0)} + z^{(k)}$. Introducing $Q_{\mathbb{K}^k}$ the least squares projection on $\mathbb{K}^k$ associated to $P^{-1}A$, $z^{(k)} = Q_{\mathbb{K}^k}\left(x^\infty - x^{(0)}\right)$ and so finally

$$x^{(k)} = Q_{\mathbb{K}^k} x^\infty + (\mathrm{Id} - Q_{\mathbb{K}^k}) x^{(0)}$$

**Proposition 8.** *GMRES is equivalent to DFPI using $Q^{(n)}$ the Least Square projection associated to $P^{-1}A$ on*

$$\mathbf{Z}^{(n)} = \left\{x^{(k)} - x^{(k-1)}\right\}_{1 \le k \le n}, \forall n \in \mathbb{N}$$

PROOF. The equality $\mathbf{Z}^{(n)} = \mathbb{K}^n$ has been shown in Prop. 6. Considering:

$$x_{\mathrm{GMRES}}^{(n)} = x^{(0)} + Q^{(n)}(x^\infty - x^{(0)})$$
$$x_{\mathrm{DFPI}}^{(n+1/2)} = x^{(n)} + Q^{(n)}(x^\infty - x^{(n)}) = x_{\mathrm{GMRES}}^{(n)} + \left(\mathrm{Id} - Q^{(n)}\right)\left(x^{(n)} - x^{(0)}\right)$$

But let's suppose that $x^{(n)} \equiv x_{\mathrm{DFPI}}^{(n)} = x^{(0)} + y^{(n)}$ with $y^{(n)} \in \mathbb{K}^n$, then

$$x^{(n+1)} = x^{(0)} + y^{(n)} + Q^{(n)}\left(x^\infty - x^{(0)} - y^{(n)}\right) + P^{-1}A\left(\mathrm{Id} - Q^{(n)}\right)\left(x^\infty - x^{(0)} - y^{(n)}\right)$$
$$= x^{(0)} + Q^{(n)}\left(x^\infty - x^{(0)}\right) + P^{-1}A\left(x^\infty - x^{(0)}\right) + P^{-1}AQ^{(n)}\left(x^\infty - x^{(0)}\right)$$
$$= x^{(0)} + y^{(n+1)}$$

where $y^{(n+1)} \in \mathbb{K}^{n+1}$. Initialization of recurrence is trivial and so $\forall n \in \mathbb{N}$:

$$\left(\mathrm{Id} - Q^{(n)}\right)\left(x^{(n)} - x^{(0)}\right) = 0$$

So finally, $x_{\mathrm{DFPI}}^{(n+1/2)} = x_{\mathrm{GMRES}}^{(n)}$. □

This means that GMRES can be viewed as a fixed-point iterations algorithm falling under the DFPI framework, which, by the way corresponds perfectly to the point of view adopted in the context of so-called Anderson acceleration method [23], whose equivalence with GMRES in the linear case has been demonstrated [24]. The framework also directly encompass deflated GMRES techniques [8].

Furthermore, as mentioned earlier BoostConv, in the linear setting, is nothing but a variant of GMRES where the least-square projection is based on $A$ instead of $P^{-1}A$.

### 4.4. Conjugate gradient methods and variants
#### 4.4.1. Conjugate gradient (CG)

CG is the Krylov method of choice for SPD matrices, as it comes with both a minimization principle and a short recurrence. As stated by Saad [2] it is nothing but: "*a realization of an orthogonal projection technique onto the Krylov subspace $\mathbb{K}^m(r_0, A)$*". In other words at each step of the algorithm, it computes the optimal solution in the norm induced by the inner product associated with $A$ within the current Krylov subspace, that is [25]:

$$x^{(k)} = \mathrm{argmin}_{y \in \mathbb{K}^k} (x^\infty - y)^T A (x^\infty - y) \tag{18}$$

The good fortune being that, contrarily to GMRES, it can be done via a short recurrence [2] thanks to the symmetry property. Expression (18) makes us realize that $x^{(k)}$ can be seen both as the Galerkin projection associated with A and the LSQ projection associated with $A^{1/2}$ of $x^\infty$ onto the Krylov subspace $\mathbb{K}^k$. Consequently, in the same manner as GMRES it falls under the DFPI framework.

*4.4.2. Bi-Conjugate gradient (BiCG)*

When A is no longer SPD however, two properties get broken wrt the Galerkin projector. First of all it is no longer a LSQ (or orthogonal) projection operator, and second of all the 3-terms relationship is lost. Trying to restore the minimization property leads to GMRES and variants, at the cost of sacrificing the 3 terms relationship; alternatively, one can decide to restore the 3-terms relationship by choosing the test subspace based on $A^T$ at the cost of loosing the minimization principle : this leads to BiCG.

More precisely in BiCG the *trial* subspace corresponds to the Krylov subspace $\mathbf{Z}^{(n)} = K^n(A, r_0)$ while the *test* subspace is chosen as $\mathbf{Y}^{(n)} = K^n(A^T, \tilde{r}_0)$ (where typically $\tilde{r}_0 = r_0$). Finally the algorithm computes at each step the Petrov-Galerkin projection of $x^\infty$, that is :

$$\text{Find } x^{(n)} \in \mathbf{Z}^{(n)}, \text{ s.t. } \forall y \in \mathbf{Y}^{(n)} \ \langle Ax^{(n)}, y \rangle = \langle Ax^\infty, y \rangle \tag{19}$$

This specific choice allows to recover the 3-terms relationship [3, 26]. Eq. (19) makes it obvious that this method falls under the DFPI umbrella with the defined Petrov-Galerkin projection. Unfortunately, as discussed in subsection 3.3, a minimization principle seems crucial to ensure smooth convergence when the projection space happens not to be invariant by the iteration matrix. That is typically the case with Krylov sequences and the observed convergence behavior of BiCG is often erratic in practice.

*4.4.3. Bi-Conjugate gradient stabilized (BiCGStab)*

This "stabilized" version of BiCG [27] is an attempt to smoothen its convergence behaviour by combining the baseline algorithm with one step of GMRES(1), that is restarted GMRES of size 1. The resulting operator is then nothing but the projection obtained by chaining BiCG Petrov-Galerkin with GMRES(1) LSQ projectors and the method enters the DFPI framework.

An improvement of convergence is indeed often observed in practice [27], although convergence is still not guaranteed [28].

## 5. Vectors Selection

Instrumental in practice for the overall algorithm's performance is the trouble vectors recruitment's strategy. Krylov methods such as BoostConv add one new vector per iteration regardless of any stability concern, while RPM happens to be more parsimonious and selective, recording Krylov vectors over a specified iteration window and then recruiting a restricted subspace out of them based on the quality of the eigenmodes approximation.

As seen earlier, Theorem 1 sheds an interesting light on this strategic difference, since RPM uses a Galerkin projection that is generally not accompanied by a minimization principle, whereas BoostConv or GMRES/Anderson acceleration use orthogonal projectors that minimize the residual in a certain norm (but not BiCG). This explains why RPM is more selective and parsimonious in its recruitment, enriching the projection base only with sufficiently well-approximated generalized eigenvectors as it requires the stability property. In contrast, BoostConv, like GMRES, indiscriminately recruits all Krylov vectors until the storage is saturated. That also explains why restarted GMRES and BoostConv lead to good results, as the dominant eigenvectors are better and better captured as iterations go by. The window size is then a limit on the maximum number of dominant eigenvalues that can be deflated.

## 6. Numerical Results

Finally we benchmark several DFPI variants on three different CFD models of relevance : a transonic wing, a transonic bump and a basic finner case. We use the HiSpAC (High-Speed Accurate and Coupled) density-based, compressible solver from the Icon-CFD softwares suite in order to perform the computations and we focus our study on the resolution of the linearized time-increment equations. In HiSpAC the system is solved monolithically (that is, in a block-coupled manner) and fully implicitly through a linear solver of choice, typically BiCGStab, which, together with standard Richardson iterations, shall constitute one of our two baselines. MILU preconditioner will be used in all compared algorithms.

Regarding the DFPI variants of interest, we will use systematically the LSQ projection associated to $A$. We will compare several trouble vectors recruitment strategies, in order to underline the theoretical results presented in the previous sections. Here is a small description of the different recruitment strategies used for the construction of the projection space:

- BoostConv : $\mathbf{Z}$ is enhanced at every iterations with the previous increment ($x^{(k-1)} - x^{(k-2)}$ at iteration $k$). No size limit. As we saw earlier, this method is similar to GMRES, down to the choice of least-squares projector for which $A$ is used instead of $P^{-1}A$

- BC-MW(x) : $\mathbf{Z}$ is enhanced at every iterations with the previous increment. A moving window of size x is used, deleting the first vector of the subspace when adding the new one if the size limit is exceeded (x chosen to store the same number of vectors as RR strategy for comparisons).

- AAOS : Add All vectors Once Stable. A temporary subspace is enhanced with the previous increment at each iteration. This subspace is not used for projection at this point. A stability test is performed on this subspace at each iteration using a Gram-Schmidt process: before adding the new vector, if it already belongs to the temporary subspace according to a given tolerance, the subspace is considered as stable/invariant. If the stability criterion is fulfilled, the entire subspace is added to the projection space and discarded, else, the new vector is added to the temporary subspace. No size limit.

- TSS : Two-Stage Stability strategy. A temporary subspace is constructed similarly to AAOS. But once stability/invariance is detected, the temporary subspace is discarded and filled in again until reaching stability again before being added to the projection space. No size limit.

- RR : Rayleigh-Ritz strategy. Again, a temporary subspace is constructed similarly to AAOS. Once stability/invariance is detected, approximate eigenvectors are calculated from this subspace using a Rayleigh-Ritz Methodology. Only well approximated eigenvectors (based on a given tolerance) are added to the projection space. No size limit.

The study will focus primarily on the iterations efficiency, that is the achieved convergence for a given number of iterations, irrespective of their individual cost. We defer the complete performance analysis to an upcoming companion paper as our main goal here is to give some insight regarding the compared convergence behaviours rather than proposing an improved algorithm.

*6.1. Transonic Wing*

The transonic wing is a steady 2D validation case using the RAE2822 profile at a free stream Mach number of 0.729, a Reynolds number of $6.5 \times 10^6$ and an angle of attack equal to $2.31°$. This is a widely studied case described in [29]. We will present results for two different pseudo-instants ($t = 33$ and $t = 66$, respectively), selected for their challenging nature for the linear solver. The CFL is ramped throughout the simulation and is worth 200 at $t = 33$ and 500 at $t = 66$. In the following figures (as well of all the subsequent ones), the Euclidean norm of the residual is plotted against the iteration number. As far as BiCGStab is concerned, the total (as opposed to the reported) number of iterations is displayed, as each step consists of a BiCG iteration followed by a GMRES(1) one.

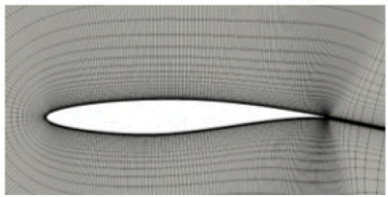

Figure 1: Transonic wing geometry and mesh.

As expected (Prop. 5), in both cases the lower and upper envelope are given, respectively, by BoostConv and Richardson methods. The AAOS strategy recovers the same residual levels as BoostConv each time the projection is updated (Prop. 7). BoostConv-MovingWindow, the Two-Stage Stability strategy, and the Rayleigh-Ritz strategy behave similarly. The convergence rate is close to the BoostConv one, although slightly lower, with significantly less storage (Table 1). TSS behaves a little better than RR as it recruits more vectors and then projects on a larger subspace but this is obviously sensitive to the selected stability threshold ($5e-2$ in our case). Finally BiCGStab convergence seems more erratic and, whilst sometimes manages to regain BoostConv convergence rate, happens to be overall slower. The salient point remains that recruiting only approximate eigenvectors maintains a convergence rate essentially similar to BoostConv (GMRES-like) while being significantly more parsimonious.

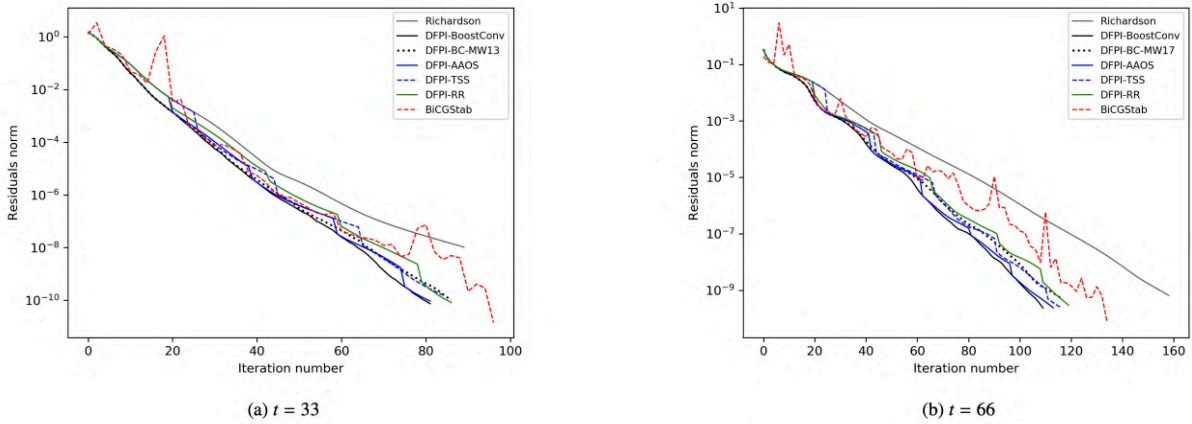

Figure 2: Transonic Wing residual comparison for different linear solvers.

Table 1: Vector storage.

|            | $t=33$ | $t=66$ |
|------------|--------|--------|
| BoostConv  | 80     | 108    |
| AAOS       | 74     | 96     |
| RR / BC-MW | 13     | 17     |
| TSS        | 23     | 22     |

## 6.2. Transonic bump

The transonic bump is a well documented validation case described in [30]. It corresponds to the simulation of the flow inside a channel with a circular bump on the lower wall. The Mach number is equal to 0.675 and the CFL to 10200.

We see Fig. 3 that the methods behave essentially in the same way as for the transonic wing, although the differences are accentuated. TSS and RR once again exhibit fast convergence, comparable to BoostConv with only a quarter of its number of vectors. BiCGStab and BoostConv-MovingWindow are further from BoostConv. Remarkably, BoostConv-MW convergence happens to be significantly degraded in this more challenging case, highlighting the importance of carefully selecting the projection vectors. Indeed RR approach selects the approximate eigenvectors whilst BoostConv-MW merely naively discards the head of the orthogonal basis.

Table 2: Vector storage.

| BoostConv  | 66 |
|------------|----|
| AAOS       | 66 |
| RR / BC-MW | 14 |
| TSS        | 16 |

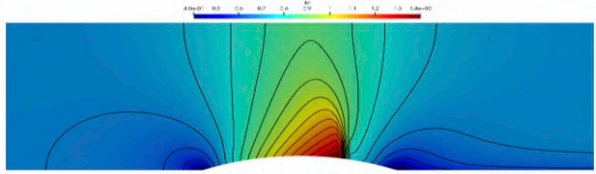
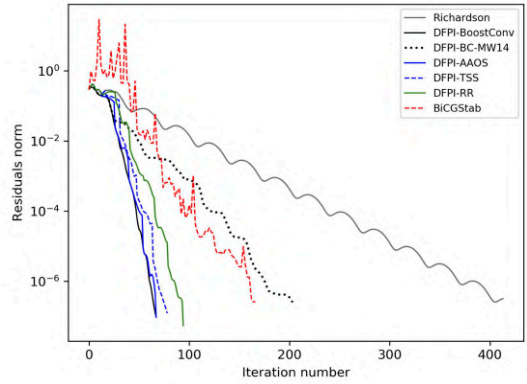

(a) Geometry and Mach contours.  (b) Residual comparison.

Figure 3: Transonic Bump.

### 6.3. Basic finner model

The last steady case corresponds to the Basic projectile with fins, named Basic Finner. It is a supersonic case described in [31] widely used for pre-industrial validation. The geometry is given below. The free stream Mach number is equal to 2.4, and the incidence is 0° at Reynolds $1.6 \times 10^6$. The CFL is equal to 500.

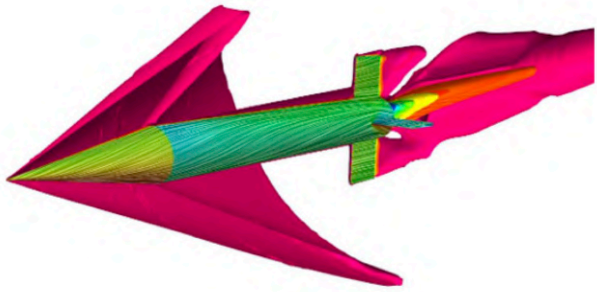
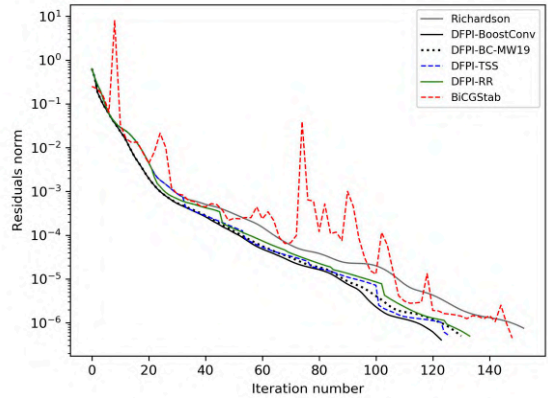

(a) Geometry and Mach contours.  (b) Residual comparison.

Figure 4: Basic finner.

Once again, the behavior of the different methods is quite similar, with an even more drastic reduction in the number of needed vectors with the TSS and RR strategies as compared to baseline BoostConv with a ratio of 1 to 6 in the case of RR!

Table 3: Vector storage.

| BoostConv | 122 |
|---|---|
| BC-MW | 19 |
| TSS | 29 |
| RR | 19 |

## 7. Conclusion

A new general framework for the iterative solution of linear systems, called DFPI, has been presented. It encompasses many existing methods such as RPM, BoostConv and popular Krylov subspaces solvers such as GMRES, PCG, BiCG-Stab and variants, highlighting the fact that they all are strikingly similar, modulo the choice of a projection operator and a projection space.

We provide a proof of convergence, irrespective of the selected operator, when the projection space is invariant by the iteration matrix or subject to small pertubations, and emphasize the importance of the availability of a minimization principle when it is not.

Finally, we discuss trouble vectors selection strategies and present numerical results on various CFD cases that support our points, a key insight being that only the invariant vectors contribute to the convergence improvement and that the remaining part of the projection space is pointless at best and detrimental in the worst case. Future work will aim at leveraging these ideas to propose an optimal DFPI algorithm delivering the best convergence rate for a given complexity overhead.

**Acknowledgments**

We would like to thank Jacques Papper and other colleagues from ICON Technology & Process Consulting, the MICS, and the FDM for the fruitful, constructive, and enlightened discussions that enabled us to progress on this project. We also thank the Academic Writing Center from CentraleSupélec for their constructive proofreading.

**Appendix A. Structure of the deflated eigenvectors**

In this section, a proof of Theorem 2 is given.

**Lemma 2.** *If the subspace $\mathbf{Z}$ is invariant under the matrix $M$, then $\mathbf{Z}$ is spanned by the initial segments of the Jordan chains of $M$. Specifically, there exist integers $q$ and $r_1, \ldots, r_q$ such that:*

$$\mathbf{Z} = \text{span}\left(\epsilon_{ij}\right)_{1 \leq i \leq q,\ 0 \leq j \leq r_i}.$$

*where $\epsilon_{ij}$ are the generalized eigenvectors of $M$.*

Proof. The matrix $M$ can be written:

$$M = \begin{pmatrix} \left(Q_R P^{-1} A Q_R\right)|_{\text{Im}(Q_R)} & B \\ 0 & D \end{pmatrix}$$

This structure arises from the invariance of $\mathbf{Z}$ under $M$, where $\mathbf{Z}$ corresponds to the image of $Q_R$.

A complete family of generalized eigenvectors of $\left(Q_R P^{-1} A Q_R\right)|_{\text{Im}(Q_R)}$ forms a basis for $\mathbf{Z}$. These eigenvectors can be extended to generalized eigenvectors of $Q_R P^{-1} A Q_R$.

Consider $e_{ij} \in \mathbf{Z}$ as generalized eigenvectors of $Q_R P^{-1} A Q_R$ associated to $\lambda_i$. We have $Q_R P^{-1} A Q_R e_{ij} = \lambda_i e_{ij} + e_{ij-1}$. Since $\mathbf{Z}$ is invariant under $P^{-1}A$, $e_{ij-1}$ must either be zero or belong to $\mathbf{Z}$. Given that $P^{-1}A\mathbf{Z} \subset \mathbf{Z}$, we have

$$P^{-1}A e_{ij} = \lambda_i e_{ij} + e_{ij-1}.$$

This shows that $e_{ij}$ is a generalized eigenvector of $M = P^{-1}A$.

By construction of the Jordan Chains and by recursion:

$$\exists \mu \in \mathbb{R}, \exists \epsilon_{ij} \in C\left(M, \lambda_i, \left(\epsilon_{ij}\right)_{0 \leq j \leq m_i - 1}\right), e_{ij} = \mu \epsilon_{ij}$$

$$\forall k \in [\![0, j-1]\!], e_{ik} = \mu \epsilon_{ik}$$

Thus, $\mathbf{Z}$ is spanned by these generalized eigenvectors: $\mathbf{Z} = \text{span}\left(\epsilon_{ij}\right)_{1 \leq i \leq q,\ 0 \leq j \leq r_i}$. □

Proof of Theorem 2. From Lemma 2, we know that $\mathbf{Z} = \text{span}\left(\epsilon_{ij}\right)_{1 \leq i \leq q,\ 0 \leq j \leq r_i}$ where $\epsilon_{ij}$ are generalized eigenvectors of $M$. Eigenvectors of $M$ that are not in $\mathbf{Z}$ are generalized eigenvectors of $M_f = (\text{Id} - Q_R) M (\text{Id} - Q_R)$ and belong to $\ker(Q_R)$. These eigenvectors come in two types:

- **Chains Starting in Z**: These are of the form $(\text{Id} - Q_R)\epsilon_{ij}$ for $1 \leq i \leq q$ and $r_i + 1 \leq j \leq m_i - 1$. The chain is initialized with $(\text{Id} - Q_R)\epsilon_{ir_i+1}$, which is an eigenvector of $M_f$. Since $\mathbf{Z}$ is invariant under $M$, we have:

$$(\text{Id} - Q_R) M (\text{Id} - Q_R)\, \epsilon_{ir_i+1} = \lambda_i\, (\text{Id} - Q_R)\, \epsilon_{ir_i+1}$$

- **Chains Not in Z**: These are complete chains $(\text{Id} - Q_R)\epsilon_{ij}$ for $q + 1 \leq i \leq n$ and $0 \leq j \leq m_i - 1$.

Consider $\epsilon_{ij} \in \mathbf{Z}$, (for $1 \leq i \leq q$ and $0 \leq j \leq r_i$):

$$N\epsilon_{ij} = Q_R \epsilon_{ij} + P^{-1} A (\text{Id} - Q_R)\, \epsilon_{ij} = \epsilon_{ij}$$

This shows that $\epsilon_{ij}$ is an eigenvector of $N$ with eigenvalue 1.

For $\epsilon_{i0}$ (where $i > q$) or $\epsilon_{ir_i+1}$ (where $i \leq q$), associated with $\lambda_i \neq 1$, define:

$$e_i^0 = (\text{Id} - Q_R)\, \epsilon_i^0 + \frac{Q_R M (\text{Id} - Q_R)\, \epsilon_i^0}{\lambda_i - 1}, \qquad N e_i^0 = \lambda_i e_i^0$$

Then, $N e_i^0 = \lambda_i e_i^0$.

The chain is constructed using:

$$\begin{cases} z_{i0} = \dfrac{Q_R M (\text{Id} - Q_R)\, \epsilon_i^0}{\lambda_i - 1} \\ e_{ij+1} = (\text{Id} - Q_R)\, \epsilon_{ij+1} + z_{ij+1} \\ z_{ij+1} = \dfrac{z_{ij} - Q_R M (\text{Id} - Q_R)\, \epsilon_{ij+1}}{1 - \lambda_i} \end{cases}$$

Thus, $N e_{ij+1} = \lambda_i e_{ij+1} + e_{ij}$.

For $\epsilon_{ij} \notin \mathbf{Z}$ (where $0 \leq j < m_i$ for $i > q$ or $r_i + 1 \leq j < m_i$ for $i \leq q$), associated with $\lambda_i = 1$, define for $2 \leq k < m_i - r_i$ where $i > q \Rightarrow r_i + 1 = 0$:

$$\begin{cases} e_{im_i-1} = (\text{Id} - Q_R)\, \epsilon_{im_i-1} + z_{im_i-1} \\ e_{im_i-k} = (\text{Id} - Q_R)\, \epsilon_{im_i-k} + Q_R M (\text{Id} - Q_R)\, \epsilon_{im_i-k+1} \end{cases}$$

For $k < m_i - r_i - 1$, we have $N e_{im_i-k} = e_{im_i-k} + e_{im_i-k-1}$ and
For $k = m_i - r_i - 1$, we have $N e_{im_i-k} = e_{im_i-k} + Q_R M (\text{Id} - Q_R)\, \epsilon_{im_i-k}$.

There are two cases:

- If $Q_R M (\text{Id} - Q_R) \epsilon_{im_i-k} = 0$, then $e_{im_i-k}$ initializes the chain.

- If $Q_R M (\text{Id} - Q_R) \epsilon_{im_i-k} \neq 0$, it is an eigenvector of $N$ in $\mathbf{Z}$, and $e_{im_i-k}$ is the second vector of a chain initialized by an eigenvector in $\mathbf{Z}$.

Thus, all eigenvectors of $N$ can be written as $e_{ij} = (\text{Id} - Q_R)\, \epsilon_{ij} + z_{ij}$, forming a free family. By dimensional arguments, there exists a bijection between the generalized eigenvectors of $N$ and those of $M$. □

**References**

[1] J.D. Anderson and J. Wendt, Computational Fluid Dynamics, Springer, 206 (1995)
[2] Y. Saad, Iterative methods for sparse linear systems, SIAM, 2003
[3] R. Fletcher, Conjugate gradient methods for indefinite systems, Numerical Analysis: Proceedings of the Dundee Conference on Numerical Analysis, 1975, Springer (2006) 73-89

[4] Y. Saad and M.H. Schultz, GMRES: A Generalized Minimal Residual Algorithm for Solving Nonsymmetric Linear Systems, SIAM Journal on Scientific and Statistical Computing, 7 (1986) 856-869

[5] B. Beckermann, A. B. J. Kuijlaars, Superlinear convergence of conjugate gradients, SIAM Journal on Numerical Analysis, 39 (2001) 300-329

[6] H.A. Van der Vorst and C. Vuik, The superlinear convergence behaviour of GMRES, Journal of Computational and Applied Mathematics, 48 (1993) 327-341

[7] R.B. Morgan, A restarted GMRES method augmented with eigenvectors, SIAM Journal on Matrix Analysis and Applications, 16 (1995) 1154-1171

[8] R.B. Morgan, GMRES with deflated restarting, SIAM Journal on Scientific Computing, 24 (2002) 20-37

[9] H. Jasak, Error analysis and estimation in the Finite Volume method with applications to fluid flows, 1996

[10] G. Karpouzas, E. Papoutsis-Kiachagias, T. Schumacher, E. De Villiers, K.C. Giannakoglou, and C. Othmer, Adjoint optimization for vehicle external aerodynamics International Journal of Automotive Engineering, 7 (2016) 1-7

[11] T. Skamagkis, E.M. Papoutsis-Kiachagias, and K.C. Giannakoglou, On the stabilization of steady continuous adjoint solvers in the presence of unsteadiness, in shape optimization, International Journal for Numerical Methods in Fluids, 93 (2021) 2677-2693

[12] P. Wesseling, Introduction to multigrid methods, Universidade Federal do Paraná, 1995

[13] R.A. Nicolaides, Deflation of conjugate gradients with applications to boundary value problems, SIAM Journal on Numerical Analysis, 24 (1987) 355-365.

[14] J. Frank and C. Vuik, On the construction of deflation-based preconditioners, SIAM Journal on Scientific Computing, 23 (2001) 442-462

[15] G.M. Shroff and H.B. Keller, Stabilization of unstable procedures: the recursive projection method, SIAM Journal on Numerical Analysis, 30 (1993) 1099-1120

[16] S. Görtz and J. Möller, Evaluation of the recursive projection method for efficient unsteady turbulent CFD simulations, 24th International Congress of the Aeronautical Sciences (2004) 1-13.

[17] V. Citro, P. Luchini, F. Giannetti, F. Auteri Efficient stabilization and acceleration of numerical simulation of fluid flows by residual recombination, Journal of computational physics, 344 (2017) 234-246

[18] A. Dicholkar, F. Zahle, N. N. Sørensen, Convergence enhancement of SIMPLE-like steady-state RANS solvers applied to airfoil and cylinder flows, Journal of Wind Engineering and Industrial Aerodynamics, 220 (2022)

[19] F.L. Bauer and C.T. Fike, Norms and exclusion theorems, Numerische Mathematik, 2 (1960) 137-141.

[20] H. Calandra, S. Gratton, J. Langou, X. Pinel, and X. Vasseur, Flexible variants of block restarted GMRES methods with application to geophysics, SIAM Journal on Scientific Computing, 34 (2012) 714-736

[21] B.W. Patton and J.P. Holloway, Application of preconditioned GMRES to the numerical solution of the neutron transport equation, Annals of Nuclear Energy, 29 (2002) 109-136

[22] S. Xu, D. Radford, M. Meyer, and J.-D. Müller, Stabilisation of discrete steady adjoint solvers, Journal of Computational Physics, 299 (2015) 175-195

[23] P. P. Pratapa, P. Suryanarayana, and J. E. Pask, Anderson acceleration of the Jacobi iterative method: An efficient alternative to Krylov methods for large, sparse linear systems, Journal of Computational Physics, 306 (2016) 43-54

[24] H.F. Walker and P. Ni, Anderson acceleration for fixed-point iterations, SIAM Journal on Numerical Analysis, 49 (2011) 1715-1735

[25] E. Paquette, T. Trogdon, Universality for the Conjugate Gradient and MINRES Algorithms on Sample Covariance Matrices, Communications on Pure and Applied Mathematics, 76 (2023) 1085-1136

[26] C. Tong, Q. Ye, Analysis of the finite precision bi-conjugate gradient algorithm for nonsymmetric linear systems, Mathematics of computation, 69 (2000) 1559-1575

[27] H. A. Van der Vorst, Bi-CGSTAB: A fast and smoothly converging variant of Bi-CG for the solution of nonsymmetric linear systems, SIAM Journal on Scientific and Statistical Computing, 13 (1992) 631-644

[28] M. H. Gutknecht, Variants of BICGSTAB for matrices with complex spectrum, SIAM journal on scientific computing, 14 (1993) 1020-1033

[29] P.H. Cook, M.A. McDonald, and M.C.P. Firmin, Aerofoil RAE2822 pressure distributions, and boundary layer and wake measurements. Experimental Data Base for Computer Program Assessment, AGARD Report AR 138, 1979

[30] R.H. Ni, A multiple-grid scheme for solving the Euler equations, AIAA Journal, 20 (1982) 1565-1571

[31] A.D. Dupuis et al., Aeroballistic range tests of the Basic Finner reference projectile at supersonic velocities, Defence Research Establishment Valcartier, (1997)